\documentclass[12pt,reqno]{amsart}
\usepackage{amsmath,amsthm}
\usepackage{amsfonts}
 \newtheorem{res}{Result}[section]

\newtheorem{theorem}[res]{Theorem}
\def\thm{\begin{theorem}}
\def\endthm{\end{theorem}}
\newtheorem{task}[res]{Problem}
\def\prob{\begin{task}}
\def\endprob{\end{task}}
\newtheorem{remark}[res]{Remark}
\def\rem{\begin{remark}}
\def\endrem{\end{remark}}
  
 \newtheorem{lemma}[res]{Lemma}
 \def\lem{\begin{lemma}}
 \def\endlem{\end{lemma}}
 
 \newtheorem{definition}[res]{Definition}
 \def\defn{\begin{definition}}
 \def\enddefn{\end{definition}}

\newtheorem{example}[res]{Example}
\def\ex{\begin{example}}
\def\endex{\end{example}}
\numberwithin{equation}{section}
\def\eqn{\begin{equation}}
\def\endeqn{\end{equation}}
\def\eqnn{\begin{eqnarray}}
\def\endeqnn{\end{eqnarray}}
\def\eqnno{\begin{eqnarray*}}
\def\endeqnno{\end{eqnarray*}}
\DeclareMathOperator*{\essinf}{essinf}
\DeclareMathOperator*{\esssup}{esssup}

\def\pf{{\sl Proof:}
\def\endpf{\phantom{.}\hfill $\square$}}

\def\ct{\Gamma}
\newcommand{\lr}[1]{\left [ #1 \right ]}
\newcommand{\lrc}[1]{\left ( #1 \right )}
\newcommand*\diff{\mathop{}\!\mathrm{d}}

\def\pay{{\mathbf v}}
\def\r{{r}}
\def\w{w}
\def\W{W}
\def\e{e}
\def\St{\mathcal{E}}

\def\al{\alpha}
\def\a{\al}

\def\M{{\mathcal M}}

\def\G{{G}}
\def\dG{{\Delta{\G}}}
\def\F{{\mathcal F}}
\def\filt{{(\Omega,\F,(\F_t)_{t\geq 0},\P)}}

\def\Bor{{\mathcal B}}

\def\R{{\mathbb R}}

\def\Om{\Omega}

\def\H {{\tilde G^*}{}}

\def\gen{{\mathcal A}}

\def\E{{\mathbb E}}
\def\P{{\mathbb P}}

\def\defto{\buildrel {\mathrm{def}}\over =}

\def\T{T}

\def\half{\frac{1}{2}}
\def\u{{\mathbf u}}
\def\v{{\mathbf v}}
\def\x{{\mathbf x}}

\def\Xm{{X^\mu}}
\def\la{{\lambda}}
\def\u{{\underline u}}
\def\v{{\underline v}}
\def\z{{\underline z}}
\def\am{{m^\alpha}}

\def\A{F}
\begin{document}
\bibliographystyle{plain}
\title[Minimising return time]{Minimising the expected commute time}
\author{Saul Jacka}
\author{Ma. Elena Hern\'andez-Hern\'andez}
\address{Department of Statistics, University of Warwick, Coventry CV4
 7AL, UK} \email{s.d.jacka@warwick.ac.uk}
 \email{m.hernandez-hernandez@warwick.ac.uk}

\begin{abstract}
Motivated in part by a problem  in simulated tempering (a form of Markov chain Monte Carlo) we seek to minimise, in a suitable sense, the time it takes a (regular) diffusion with instantaneous reflection at 0 and 1  to travel
from the origin to $1$ and then return (the so-called commute time from 0 to 1). 
 We consider the static and dynamic versions of this problem where the control mechanism is related to the diffusion\rq{}s drift via the corresponding scale function. In the static version the  diffusion\rq{}s drift can be chosen at each point in [0,1], whereas in the dynamic version, we are only able to choose the drift at each point at the time of first visiting that point. The dynamic version leads to a novel type of stochastic control problem.

\end{abstract}

\thanks{ Saul Jacka gratefully acknowledges funding received from the EPSRC grant EP/P00377X/1 and is also grateful to the Alan Turing Institute for their financial support under the EPSRC 
grant EP/N510129/1.}
\thanks{The authors are most grateful to Gareth Roberts for suggesting this problem}

\thanks{ The authors would like to thank David Hobson, Jon Warren, Martin Barlow and Sigurd Assing for stimulating discussions on this problem}
\maketitle \centerline{{\today}}

\section{Introduction and problem motivation}
\subsection{Introduction} 

Suppose that $X^\mu$ is the diffusion on $[0,1]$  started at 0 and given by
\begin{equation}\label{SDE}
\diff X^\mu_t=\sigma(X^\mu_t) \diff B_t+\mu(X^\mu_t) \diff t\,\,\,\text{ on $\,$ (0,1)}
\end{equation}
with instantaneous reflection at 0 and 1 (see \cite{R+W} or \cite{IM} for details).
 Where there is no risk of confusion we omit the superscript $\mu$.

Formally, we define $\T_x$ to be the first time that the diffusion
reaches $x$, then we define $\ct=\ct(X)$, the commute time (between 0 and 1),
by
$$
\ct(X)\defto \inf\left \{ t>\T_1(X):\: X_t=0 \right \}.
$$
The commute time is defined  for random walks on graphs in \cite{Barlow}. The original commute time identity (which we give later in (\ref{comm})) was only discovered in 1989 and first appeared in \cite{CRRST}

In this article, we consider the following problem (and several variants and generalisations):
\prob\label{prob1} Minimise the expected commute time $\E [\, \ct\,]$; i.e. find
$$
\inf_\mu\E [\ct(X^\mu)],
$$
where the infimum is taken over a suitably large class of drifts $\mu$, to be specified in more detail later.
\endprob
Given the symmetry of the problem it is tempting to hypothesise that the optimal choice of $\mu$ is 0. We shall soon see that, in general, this is false, although it {\em is} true when $\sigma\equiv 1$

We actually want to try and minimise $\ct$ (and additive functionals of $\Xm$ evaluated at $\ct$) in as general a way as possible so we extend Problem \ref{prob1} in the following ways:
\prob\label{prob2}
Find
$$
\inf_\mu \E\lr{ \int_0^\ct f(X^\mu_t) \diff t}
$$
for suitable positive functions $f$;
\endprob
and
\prob\label{prob3}
Find
$$
\sup_\mu \E \lr{\exp \lrc{-\int_0^\ct\alpha(X^\mu_t) \diff t}};
$$
for suitable positive functions $\a$.
\endprob
Although we will prove more general versions it seems appropriate to give a preliminary statement of results in this context.
\thm\label{thm1}
Suppose that $\sigma$ is a strictly positive function on $[0,1]$,  that $f$ is a non-negative Borel-measurable function on $[0,1]$ and that, denoting Lebesgue measure by $\la$ ,
\eqn\label{ass1}
\frac{\sqrt{f}}{\sigma}\in L^1([0,1],\la), 
\endeqn
 then 
$$
\inf_{\text{measurable }\mu}\E \lr{\int_0^\ct f(\Xm_t) \diff t}=\lrc{\int_0^1\sqrt{\frac{2f(u)}{\sigma^2(u)} }\diff u}^2.
$$
If, in addition, $\sqrt{\frac{f}{\sigma^2} }$ is continuously differentiable and strictly positive on (0,1), then the optimal drift is $\hat \mu$ given by
$$
\hat \mu=-\half \lrc{\ln \sqrt{\frac{f}{\sigma^2} }}\rq{}.
$$
\endthm

\thm\label{thm2}
Suppose that $\sigma$ is  a strictly positive function on $[0,1]$,  that $\a$ is a non-negative Borel-measurable function on $[0,1]$ and that 
\eqn\label{ass2}
\frac{\sqrt{\a}}{\sigma}\in L^1([0,1],\la),
\endeqn
$$
\sup_{\text{measurable }\mu}\E \lr{\exp\lrc{-\int_0^S\a (\Xm_t) \diff t}}=\cosh\lrc{\int_0^1\sqrt{ \frac{2\a(u)}{\sigma^2(u)}} \diff u}^{-2}.
$$
If, in addition, $\sqrt{ \frac{\a}{\sigma^2} }$ is continuously differentiable and strictly positive on (0,1), then the optimal drift is $\hat \mu$ given by
$$
\hat \mu=-\half \lrc{\ln \sqrt{ \frac{\a}{\sigma^2} }}\rq{}.
$$
\endthm

We will eventually solve the problems dynamically, i.e. we will solve the corresponding stochastic control problems. However, we shall need to be careful about what these are as the problem is essentially non-Markovian. 
Normally in stochastic control problems, one can choose the drift of a controlled diffusion at each time point (see, e.g., \cite{FSoner}) but this is not appropriate here. In this context, it is appropriate that the drift is \lq frozen\rq{} once we have had to choose it for the first time. We shall formally model this in section \ref{stoch} .

\subsection{Problem motivation}
The problem gives a much simplified model for one arising in {\em simulated tempering} -a form of Markov Chain Monte Carlo (MCMC). Essentially the level corresponds to the temperature in a \lq\lq{}heat bath\rq\rq{}. 
The idea is that when simulating a draw from a highly multimodal distribution we use a reversible Markov Process to move between low and high  temperature states (and thus smear out the modes temporarily) so that the Markov chain can then move around the statespace; then at low temperature we sample from the true distribution (see \cite{ARR}).

The rest of the paper is organised as follows. The next section introduces some notation and preliminary results. Section \ref{StaticSection} contains the static (generalised) versions of Problem \ref{prob2} and Problem \ref{prob3}. The dynamic versions are presented in Section  \ref{stoch}. Then, in Section \ref{DiscrC},  we solve the corresponding discrete statespace problems (in both discrete and continuous time). Some examples are given in Section \ref{S-Ex}. We provide the proofs of the main results in an appendix.

\section{Notation and some general formulae}
We need to define the set of admissible controls quite carefully and two approaches suggest themselves: the first is to restrict controls to choosing the drift $\mu$  whilst the second is to control the corresponding random scale function. In the interest of generality, we adopt the second approach.

We assume the usual Markovian setup, so that each stochastic process lives on a suitable filtered space $(\Om,\F,(\F_t)_{t\geq 0})$, with the usual family of probability measures $(\P_x)_{x\in[0,1]}$ corresponding to the possible initial values.

Let $X^{\mu}$ be the diffusion with instantaneous reflection as given in \eqref{SDE}. Denote by $s^\mu$ the standardised \textit{scale function} of   $X^\mu$  and by $m^\mu$ the corresponding \textit{speed measure}. 
\begin{remark}
Since $X^\mu$ is regular and reflection is instantaneous we have:
$$
s^\mu(x) = \int_0^x\exp \lrc{-2\int_0^u\frac{\mu(t)}{\sigma^2(t)} \diff t} \diff u,
$$
$$
m^\mu([0,x])\defto m^\mu(x)=2\int_0^x\frac{ \diff u}{\sigma^2(u)s\rq{}(u)}=2\int_0^x\frac{\exp \lrc{2\int_0^u\frac{\mu(t)}{\sigma^2(t)} \diff t}}{\sigma^2(u)} \diff u,
$$
(see \cite{RY}).
\end{remark}
From now on, we shall consider the more general case where we only know that (dropping the $\mu$ dependence) $s$ and $m$ are absolutely continuous with respect to Lebesgue measure so that, denoting the respective Radon-Nikodym derivatives by $s\rq{}$ and $m\rq{}$ we have
$$
s\rq{}m\rq{}= \frac{2}{\sigma^2}\,\,\text{ $\lambda$-a.e.}
$$
For such a pair we shall denote the corresponding diffusion by $X^s$. We underline that we are only considering regular diffusions with \lq\lq{}martingale  part\rq\rq{} $\int\sigma \diff B$ or, more precisely, diffusions $X$ with scale functions $s$ such that
$$
 \diff s(X_t)=s\rq{}(X_t)\sigma(X_t) \diff B_t,
$$
so that, for example, sticky points are excluded (see \cite{E} for a description of driftless sticky BM and its construction, see also \cite{E2} for other problems arising in solving stochastic differential equations ). 
\begin{remark}
Note that our assumptions do allow generalised drift: if $s$ is the difference between two convex functions (which we will not necessarily assume) then 
\eqn\label{gend}
X^s_t=x+\int_0^t\sigma(X^s_u) \diff B_u-\half\int_{\R}L^a_t(X)\frac{s\rq{}\rq{}( \diff a)}{s\rq{}_-(a)},
\endeqn
where $s\rq{}_-$ denotes the left-hand derivative of $s$,  $s\rq{}\rq{}$ denotes the signed measure induced by $s\rq{}_-$ and $L^a_t(X)$ denotes the local time at $a$ developed by time $t$ by $X$ (see \cite{RY} Chapter VI for details).
\end{remark}
\defn
For each $y\in [0,1]$, we denote by $\phi_y$ the function
$$
\phi_y:\; x\mapsto \E_x \lr{\int_0^{\T_y}f(X_t)\diff t},
$$
where, as is usual, the subscript $x$ denotes the initial value of $X$ under the corresponding law $\P_x$.
\enddefn
\thm\label{time}
For $0\leq x\leq y$, $\phi_y$ is given by
\eqn\label{time1}
\phi_y(x)=\int_x^y\int_{u=0}^vf(u)m\rq{}(u)s\rq{}(v)\diff u \diff v,
\endeqn
while for $0\leq y\leq x$, $\phi_y$ is given by
\eqn\label{time2}
\phi_y(x)=\int_y^x\int_{u=v}^1f(u)m\rq{}(u)s\rq{}(v) \diff u \diff v.
\endeqn
In particular,
\eqn\label{comm}
\E_0\lr{\int_0^\ct f(X^s_t) \diff t}=\int_0^1\int_0^1f(u)m\rq{}(u)s\rq{}(v) \diff u \diff v.
\endeqn
\endthm
\pf
This follows immediately from 
Proposition VII.3.10 of \cite{RY} on observing that, with instantaneous reflection at the boundaries, the speed measure is continuous.
\endpf

We give similar formulae for the discounted problem:
\defn
We denote by $\psi_y$ the function
$$
\psi_y:\; x\mapsto \E_x\lr{\exp\lrc{ -\int_0^{\T_y}\a (X_t) \diff t }}.
$$
\enddefn
For each $n$, denote by $I_n(x)$ the integral
$$
I_n(x)\,\,\defto\int\limits_{0\leq u_1\leq v_1\leq u_2\ldots \leq v_n\leq x}\a(u_1)\ldots \a(u_n)\diff m(u_1)\ldots \diff m(u_n)\diff s(v_1)\ldots \diff s(v_n),
$$
and by $\tilde I_n(x)$  the integral
$$
\tilde I_n(x)\,\,\defto\int\limits_{x\leq v_1\leq u_1\leq v_2\ldots \leq u_n\leq 1}\a(u_1)\ldots \a(u_n) \diff m(u_1)\ldots \diff m(u_n)\diff s(v_1)\ldots \diff s(v_n),
$$
with $I_0=\tilde I_0\equiv 1$.
Define $G$ and $\tilde G$ by
\eqn\label{disc}
G(x)\defto \sum_0^\infty I_n(x) \,\,\,\text{ and }\,\,\,\tilde G(x) \defto \sum_0^\infty \tilde I_n(x).
\endeqn
\thm\label{disc1}
\begin{itemize}
\item[(i)]
Either
\eqn\label{bound}
\int_0^1 \a(u) \diff m(u)<\infty, 
\endeqn
or
$$\E_0 \lr{\exp\lrc{ -\int_0^{\ct}\a (X_t) \diff t }}=0,$$
in which case 
$$
\int_0^{\ct}\a (X_t) \diff t=\infty \text{ a.s.}
$$
\item[(ii)]
Now suppose that (\ref{bound}) holds. 
 Then, the sums in (\ref{disc}) are convergent and for $x\leq y$
$$
\psi_y(x)=\frac{G(x)}{G(y)},
$$
while for $x\geq y$
$$
\psi_y(x)=\frac{\tilde G(x)}{\tilde G(y)}.
$$
\end{itemize}
\endthm
\goodbreak

\section{The static control problems}\label{StaticSection}
For now we will state and prove more general, but still {\em non-dynamic versions} of Theorems \ref{thm1} and \ref{thm2}.

We define our constrained control set as follows:
\defn
Given a fixed scale function $s_0\sim \la$ and $C$, a Borel subset of [0,1], we define the constrained control set $\M^C_{s_0}$ by
\eqn\label{constr}
\M^C_{s_0}=\left \{\text{ scale functions }\, s\,:\; \diff s|_C= \diff s_0|_C\,\text{ and } s\sim \la \right \}.
\endeqn
For each $s \in \M^C_{s_0}$, the corresponding controlled diffusion $X^s$ has scale function $s$ and speed measure $m$ given by 
$$
m\rq{}=\frac{2}{\sigma^2 s\rq{}}.
$$
\enddefn
\thm\label{thm3}
For any scale function $s\sim\la$, define the measure $I^s$ on $([0,1],\Bor([0,1])$ by
$$
I^s(D)\defto \int_Df(u)m(\diff u)=\int_D 2\frac{f(u)}{\sigma^2(u)s\rq{}(u)} \diff u
$$
and the measure $J$ by
$$
J(D)\defto\int_D\sqrt{\frac{2f(u)}{\sigma^2(u)}} \diff u,
$$
then, given a scale function $s_0$,
$$
\inf_{s\in \M^C_{s_0}  }\E_0\lr{\int_0^\ct f(X^s_t)\diff t}=\biggl(\sqrt{s_0(C)I^{s_0}(C)}+J(C^c)\biggr)^2.
$$

The optimal choice of $s$ is given by 
$$
s(\diff x)=\begin{cases}
s_0(\diff x)&:\text{ on }C\\
\sqrt{\frac{s_0(C)}{I^{s_0}(C)}}\sqrt{\frac{2f(x)}{\sigma^2(x)}} \diff x&: \text{ on } C^c
\end{cases}$$

\endthm

For the generalised version of  the discounted case (Problem \ref{prob3}) we only deal with constraints on $s$ on $[0,y]$.

\thm\label{discstat}
Assume that (\ref{bound}) holds and
define 
$$
\tilde\sigma^2(x)\defto \frac{\sigma^2(x)}{\a (x)},
$$
\begin{itemize}
\item[(i)]Let $G$ be as in equation (\ref{add3}), so that (at least formally)
$$
\half \alpha\lrc{\tilde\sigma^2s\rq{}\lr{\frac{G\rq{}}{s\rq{}}}\rq{} -2G }=\half \sigma^2 G\rq{}\rq{}+\mu G\rq{}-\alpha G=0
$$
and let $\H $ satisfy the \lq\lq{}adjoint equation\rq\rq{}
$$\half\alpha\bigl( \frac{\bigl[\tilde\sigma^2s\rq{}\H \rq{}\bigr]\rq{}}{s\rq{} }-2\H \bigr)=0
$$
with boundary conditions $\H (0)=1\hbox{ and }\H \rq{}(0)=0$, so that
\eqnn\label{H}
\H (x)&=&1+\int_{v=0}^x\int_{u=0}^v\frac{2\a(v)\H (v)}{\sigma^2(v)s\rq{}(v)}s\rq{}(u)\diff u \diff v\\
&=&1+\int_{v=0}^x\int_{u=0}^v \a(v)\H(v) \diff m(v) \diff s(u),\nonumber
\endeqnn
then $\H $ is given by
\eqn
\H (x)=\sum_{n=0}^\infty \tilde I_n^*(x), 
\endeqn
where
\eqn
\tilde I^*_n(x)\defto\int\limits_{0\leq u_1\leq v_1\leq\ldots v_n\leq x}
\a(v_1)\ldots \a(v_n)\diff s(u_1)\ldots \diff s(u_n)\diff m(v_1)\ldots \diff m(v_n).
\endeqn
\item[(ii)]The optimal payoff for Problem \ref{prob3} is given  by
$$
\sup_{s\in \M^{[0,y]}_{s_0} }\E_0\lr{\exp\lrc{-\int_0^\ct \a(X^s_t)dt}}=\hat \psi(y)
$$
where
$$\hat\psi(y)=\lrc{\sqrt{G\H }\cosh \A(y)+\sqrt{\tilde\sigma^2G\rq{}\H \rq{}}\sinh \A(y)}^{-2},$$
with
$$
\A(y)\defto \int_y^1 \frac{\sqrt{2}\diff u }{\tilde\sigma(u)}=\int_y^1 \sqrt{ \frac{2\a}{\sigma^2(u)} } \diff u.
$$
The payoff is attained by setting $\tilde\sigma(x)s\rq{}(x)=\sqrt{\frac{G\H \rq{}}{G\rq{}\H }(y)}$ for all $x\geq y > 0$.  If $y=0$, any {\em constant} value for   $\tilde\sigma(x)s\rq{}(x)$ will do.\\
\end{itemize}
\endthm

\rem\label{skew}
We see that, in general, in both Theorems \ref{thm3} and \ref{discstat} the optimal scale function has a discontinuous derivative.
 In the case where $C=[0,y)$ there is a discontinuity in $s\rq{}$ at $y$. This will correspond to partial reflection at $y$ (as in skew Brownian motion- see \cite{RY} or \cite{IM}) and will give rise to a singular drift -- at least at $y$.
\endrem
\rem
We may easily extend Theorems \ref{thm3} and \ref{discstat} to the cases where  $f$ or $\a$ vanishes on some of [0,1]. 
In the case where $N\defto\{ x:\; f(x)=0\}$ is non-empty, observe first that 
the cost functional does not depend on the amount of time the diffusion spends in $N$ so that every value for $\diff s|_N$ which leaves the 
diffusion recurrent will give the same expected cost. If $\la(N)=1$ then the problem is trivial, otherwise  define the revised statespace $\mathcal{E}=[0,1-\la(N)]$ 
and solve the problem on this revised interval with the cost function $\tilde f(x)\defto f(g^{-1}(x))$ where
$$
g: t\mapsto \la([0,t]\cap N^c)
$$
and
$$
g^{-1}:x\mapsto\inf\{ t:\;g(t)=x \}.
$$
This gives us a diffusion and scale function $s^\mathcal{E}$ which minimises the cost functional on $\mathcal{E}$. Then we can extend this to a solution of the original problem  by taking
$$
\diff s= \diff s_01_N+ \diff \tilde s1_{N^c},
$$
where
$ \diff s_0$ is any finite measure equivalent to $\la$ and 
$ \diff \tilde s$ is the Lebesgue-Stiltjes measure given by
$$
\tilde s([0,t])=s^\mathcal{E} ([0,\la([0,t]\cap N^c)=s^\mathcal{E} (g(t)).
$$

An exactly analogous method will work in the discounted problem
\endrem

\section{The dynamic control problems}\label{stoch}
We now turn to the dynamic (generalised) versions of Problems \ref{prob2} and \ref{prob3}.

A  moment\rq{}s consideration shows that it is not appropriate to model the dynamic version of the problem by allowing the drift to be chosen adaptively.
If we were permitted to do this then we could choose a very large positive drift until the diffusion reaches 1 and then a very large negative 
drift to force it back down to 0. The corresponding optimal payoffs for Problems  \ref{prob2} and \ref{prob3} would be 0 and 1 respectively. 
We choose, instead, to consider the problem where the drift may be chosen dynamically at each level, but only  when the diffusion first reaches 
that level. Formally, reverting to the finite drift setup, we are allowed to choose controls from the collection $\M$ of adapted processes $\mu$ 
with the constraint that 
\eqn\label{drift}
\mu_t=\mu_{\T_{X_t}},
\endeqn
 or continuing the generalised setup, to choose scale measures dynamically, in such a way that $s\rq{}(X_t)$ is adapted.

Although these are very non-standard control problems we are able to solve them -- mainly because we can exhibit an explicit solution-- following the same control as in the \lq\lq{}static\rq\rq{} case. 
\rem
Note that this last statement would not be true if our constraint was not on the set $[0,y]$. 
To see this, consider the case where our constraint is on the set $[y,1]$. If the controlled diffusion starts at $x>0$ then there is a positive 
probability that it will not reach zero before hitting 1, in which case the drift will not have been chosen at levels below $I_{\T_1}$, 
the infimum on $[0, \T_1]$. Consequently, on the way down we can set the drift to be very large and negative below $I_{\T_1}$
Thus the optimal dynamic control will achieve a strictly lower payoff than the optimal static one in this case.
We do not pursue this problem further here but intend to do so in a sequel.
\endrem

As pointed out before,  two approaches for the set of admissible controls are available: the first is to restrict controls to choosing a drift with the property (\ref{drift}) whilst the second is to allow suitable random scale functions.

Both approaches have their drawbacks: in the first case we know from Remark \ref{skew} that, in general, the optimal control will not be in this class, whilst, in the second, it is not clear how large a class of random scale functions will be appropriate. In the interests of ensuring that an optimal control exists, we again adopt the second approach. From now on, we fix the Brownian Motion $B$ on the filtered probability space $\filt$. 
\defn
By an {\em equivalent random scale function} we mean a random measure $s$ belonging to the class $\M$ defined by
\eqnn
&\M\defto &\Big \{\text{random, finite Borel measure $s$ on $[0,1]$  }\,:\,\,\, s\sim \la \text{ and } \Big. \\
&& \Big .\phantom{m}\text{there exists a {\em martingale} }Y^s\text{ with }Y^s_t=\int_0^ts\rq{}\circ s^{-1}.\sigma\circ s^{-1}(Y_u) \diff B_u\Big \}.\nonumber
\endeqnn
For each $s\in \mathcal{M}$, we define the corresponding controlled process $X^s$ by
$$
X^s_t=s^{-1}(Y_t).
$$
\enddefn
\rem In general, the martingale constraint is both about existence of a solution to the corresponding stochastic differential equation and about imposing a suitable progressive measurability condition on the random scale function.
\endrem
We define our constrained control set as follows:
\defn
Given a fixed scale function $s_0\in \mathcal{M}$ and $C$, a Borel subset of [0,1], we define the constrained control set $\M^C_{s_0}$ by
\eqn\label{constr}
\M^C_{s_0}=\left \{ s\in \mathcal{M}:\; \diff s|_C= \diff s_0|_C \right \}
\endeqn
\enddefn
\rem
Note that $\M$ contains all {\em deterministic} equivalent scale functions. An example of a random element of $\M$ when $\sigma\equiv 1$ is $s$, given by
$$
ds|_{[0,\half)}=d\la;\; \diff s(x)|_{[\half, 1]} =1\diff \la 1_{(\T_{\half}(B)<1)}+\exp(-2(x-\half))d\la 1_{(\T_{\half}(B)\geq 1)},
$$
corresponding to $X^s$ having drift 1 above level $\half$ if and only if it reaches that level before time 1.

\endrem

\thm\label{dsc1}
For each $s\in\M$, let $M^s_t$ denote the running maximum of the controlled process $X^s$.
Then for each $s_0\in \M$, the optimal payoff (or Bellman) process $
V^{s_0}$ defined by
$$
V^{s_0}_t\defto \essinf\limits_{ s\in \M^{s_0}_{M^{s_0}_t} }\E\lr{\int_0^\ct f(X^s_t)\diff t \Big |\F_t}
$$
is given by
\eqn\label{soln1}
V^{s_0}_t=v_t\defto \begin{cases}
\int_0^tf(X^{s_0}_u)\diff u+2\bigl(\sqrt{{s_0}(M^{s_0}_t)I^{s_0}(M^{s_0}_t)}+J(M^{s_0}_t)\bigr)^2-\phi_{X^{s_0}_t}(0):&\text{ for }M^{s_0}_t<1\\
\int_0^{t\wedge \ct}f(X^{s_0}_u)\diff u+\phi_0(X^{s_0}_{t\wedge \ct}):&\text{ for }M^{s_0}_t=1,
\end{cases}
\endeqn
(where $\phi$ is formally given by equations (\ref{time1}) and (\ref{time2}) with $s=s_0$),
and the optimal control is to take 
\eqn\label{soln2}
s\rq{}(x)=\frac{s_0(M^{s_0}_t)}{I^{s_0}(M^{s_0}_t)}\frac{\sqrt{f(x)}}{\sigma(x)}\,\,\, \text{ for }\,\, x\geq M^{s_0}_t.
\endeqn
\endthm

\thm\label{dsc2}
The Bellman process for Problem \ref{prob3} is given  by
$$V^{s_0}_t\defto \esssup\limits_{ s\in \M^{s_0}_{M^{s_0}_t} }\E\lr{\exp\lrc{-\int_0^\ct \a(X^{s_0}_t) \diff t}\Big | \F_t }=\pay_t\defto e^{-\int_0^t \alpha(X^{s_0}_u)\diff u}\psi(X^{s_0}_t,M^{s_0}_t),$$
where
\begin{equation*}\psi(x,y)=\begin{cases}G(x)\hat\psi(y)&\text{if }y<1,\\
\tilde G(x)&\text{if }y=1.
\end{cases}
\end{equation*}
and
$$\hat\psi(y)=\lrc{\sqrt{G\H }\cosh \A(y)+\sqrt{\tilde\sigma^2G\rq{}\H \rq{}}\sinh \A(y)}^{-2},$$
(as in Theorem \ref{discstat} (ii)) with
$$
\A(y)\defto \int_y^1 \frac{\diff u}{\tilde\sigma(u)}.
$$
The payoff is attained by setting 
\eqn\label{d3}
\tilde\sigma(x)s\rq{}(x)=\sqrt{\frac{G\H \rq{}}{G\rq{}\H }(M^{s_0}_t)}\,\,\,\text{ for all }\,\,x\geq M^{s_0}_t.
\endeqn
\endthm

\section{The discrete statespace case}\label{DiscrC}
\subsection{Additive functional case}
Suppose that $X$ is a discrete-time birth and death process on $E=\{0,\ldots,N\}$, with transition matrix $(P)$ given by
$$
P_{n,n+1}=p_n\text{ and }1-p_n=q_n=P_{n,n-1}\text{ and }p_N=q_0=0.
$$
We define
$$\w_n\defto\frac{q_n}{p_n}
$$
and
$$\W_n=\prod_{k=1}^n\w_k,
$$
with the usual convention that the empty product is 1.
Note that $s$, given by 
$$
s(n)\defto \sum_{k=0}^{n-1} \W_k,
$$
is the discrete scale function in that $s(0)=0$, $s$ is strictly increasing on $E$ and 
$s(X_t)$ is a martingale.
\rem Note that when we choose $p_n$ or $\w_n$ we are implicitly specifying $s(n+1)-s(n)$ so we shall denote this quantity by $\Delta s(n)$ and we shall denote by $\diff s$ the Lebesgue-Stiltjes measure on $\St \defto\{0,1,\ldots,N-1\}$ given by
$$
\diff s(x)=\Delta s(x).
$$
\endrem

Let $f$ be  a positive function on $E$ and define
$$
\tilde f(n)=\half(f(n)+f(n+1))\text{ for }0\leq n\leq N-1.
$$
\thm
If we define
$$
\phi_y(x)=\E_x \lr{\sum_{t=0}^{\T_y-1}f(X^{s_0}_t)},
$$
then
for $x\leq y$
\eqn\label{time3}
\phi_y(x)=f(x)+\ldots+f(y-1)+\sum_{v=x}^{y-1}\sum_{u=0}^{v-1}\frac{2\tilde f(u)W_v}{W_u}
\endeqn
while for $y\leq x$
\eqn\label{time4}
\phi_y(x)=f(y+1)+\ldots+f(x)+\sum_{v=y}^{x-1}\sum_{u=v+1}^{N-1}\frac{2\tilde f(u)W_v}{W_u}.
\endeqn
\endthm
\pf
It is relatively easy to check that $\phi$
satisfies the linear recurrence
$$
\phi(x)=p_x\phi(x+1)+q_x\phi(x-1)+f(x)
$$
with the right boundary conditions.
\endpf

It follows from this that
optimal payoffs are given by essentially the same formulae as in the continuous case.
Thus we now define
 the constrained control set $\M^C_{s_0}$ by
\eqn\label{constr2}
\M^C_{s_0}=\{\text{ scale functions }s:\; \diff s|_C=\diff s_0|_C \,\text\}.
\endeqn
\rem By convention we shall always assume that $0\in C$ since we cannot control $\W_0$ and hence cannot control $\Delta s(0)$.
\endrem
\thm\label{discrete1}
For any scale function $s$, define the measure $I^s$ by
$$
I^s(D)\defto \sum_{k\in D}\frac{2\tilde f(k)}{W_k}\text{ for }D\subseteq\St
$$
and the measure $J$ by
$$
J(D)\defto\sum_{k\in D}\sqrt{2\tilde f(k)}\text{ for }D\subseteq\St,
$$
then, given a scale function $s_0$,
$$
\inf_{s\in \M^C_{s_0}  }\E_0\lr{\sum_0^\ct f(X^s_t)\diff t}=\lrc{\sqrt{s_0(C)I^{s_0}(C)}+J(C^c)}^2.
$$
The optimal choice of $s$ is given by 
$$
\Delta s(x)=\W_x=\begin{cases}
\W^0_x&:\text{ on }C\\
\sqrt{\frac{s_0(C)}{I^{s_0}(C)}}\sqrt{2\tilde f(x)}&: \text{on } C^c
\end{cases}$$
\endthm
\rem Note that all complements are taken with respect to $\St$.
\endrem

The dynamic problem translates in exactly the same way:
we define the constrained control set:
$$
\M^{s_0}_y=\{ s:\; \diff s|_{\{0,\ldots, y-1\}}=\diff s_0|_{\{0,\ldots, y-1\}}\},\: y\geq 1
$$
then we have the following:
\thm\label{discrete2}
For each $s\in\M$, let $M^s_t$ denote the running maximum of the controlled process $X^s$.
Then for each $s_0\in \M$, the optimal payoff (or Bellman) process $
V^{s_0}$ defined by
$$
V^{s_0}_t\defto \essinf_{ s\in \M^{s_0}_{M^{s_0}_t} }\E\lr{\sum_0^{\ct-1}f(X^s_t)\diff t\, \Big |\, \F_t}
$$
is given by
\eqn\label{soln3}
V^{s_0}_t=v_t\defto \begin{cases}
\sum\limits_0^tf(X^{s_0}_u)+2\lrc{\sqrt{{s_0}(M^{s_0}_t)I^{s_0}(M^{s_0}_t)}+J(M^{s_0}_t)}^2-\phi_{X^{s_0}_t}(0):&\text{ for }M^{s_0}_t<N\\
\sum\limits_0^tf(X_u)du+\phi_0(X^{s_0}_t):&\text{ for }M^{s_0}_t=N,
\end{cases}
\endeqn
(where $\phi$ is formally given by equations (\ref{time3}) and (\ref{time4}) with $s=s_0$),
and the optimal control is to take 
\eqn\label{soln4}
\W(x)=\frac{s_0(M^{s_0}_t)}{I^{s_0}(M^{s_0}_t)}\sqrt{\tilde f(x)}\,\,\, \text{ for }\,\, x\geq M^{s_0}_t.
\endeqn
\endthm

\subsection{The discounted problem}
Suppose that for each $x\in \St$, $0\leq \r_x\leq 1$, then 
define
$$\sigma^2_i\defto (1-\r_{i-1}\r_i)^{-1}, \hbox{ with }\r_{-1}\hbox{ taken to be }1,$$
and
$\sigma(i_1,i_2,\ldots, i_l)\defto \prod\limits_{m=1}^l \sigma_{i_m}$.
Now set 
$$A_k(x)=\{(\u,\v):\; 0\leq u_1<v_1<\ldots <v_k<x\},
$$
$$\tilde A_k(x)=\{(\u,\v):\; x\leq v_1<u_1<\ldots <u_k<N\}
$$ 
and $\W^\sigma_m=\sigma_m\W_m$, where
$\W_m$ is as before. 
Note that $A_k$ and $\tilde A_k$ will be empty for large values of $k$.
\begin{theorem}\label{pgf}
For $x\leq y$
\eqn\label{h3}
E_x\lr{\prod_{t=0}^{T_y-1}\r_{X^{s_0}_t}}=\r_x\ldots \r_{y-1}\G(x)/ \G(y),
\endeqn
where 
$$
\G(x)\defto 
1+\sum_{k=1}^\infty \sum_{(\u,\v)\in A_k(x)}\frac{1}{\sigma(\u,\v)}\prod_{m=1}^k\frac{\W^\sigma_{v_m}} {\W^\sigma_{u_m}},
$$
while for $x\geq y$
\eqn\label{h4}
E_x \lr{\prod_{t=0}^{T_y-1}\r_{X^{s_0}_t}}=\r_{y+1}\ldots \r_{x}\tilde \G(x)/ \tilde \G(y)
,
\endeqn
where
$$\tilde \G(x)\defto 
1+\sum_{k=1}^\infty \sum_{(\u,\v)\in \tilde A_k(x)}\frac{1}{\sigma(\u,\v)}\prod_{m=1}^k\frac{\W^\sigma_{v_m}} {\W^\sigma_{u_m}},
$$
\end{theorem}

\pf
Define $d_x=\E_x\prod\limits_{t=0}^{\T_{x+1}-1}r(X^{s_0}_t)$, then
$$
d_x=\r_x(p_x+q_xd_{x-1}d_x).
$$
Setting
$$
\r_xt_x/t_{x+1}=d_x
$$
we see that
$$
t_{x+1}=(1+w_x)t_x-w_x\r_x\r_{x-1}t_{x-1}
$$
or
$$
(t_{x+1}-t_x-w_x(t_x-t_{x-1}=w_x(1-\r_{x-1}\r_x)t_{x-1}.
$$
Substituting
$$t_x=\G(x)
$$
It is easy to check that this is satisfied.
Now boundary conditions give equation (\ref{h3}).
The proof of equation (\ref{h4}) is essentially the same.
\endpf

Now with this choice of $\sigma$
we get same results as before:

\thm\label{discrete3}
Suppose $\G$ and $\tilde \G$ are as defined in Theorem \ref{pgf}:
we set
$$\G^*(x)\defto 1+\sum_{k=1}^\infty \sum_{(\u,\v)\in A^*_k(x)}\frac{1}{\sigma(\u,\v)}\prod_{m=1}^k\frac{\W^\sigma_{v_m}} {\W^\sigma_{u_m}}
$$
and
$$\tilde \G^*(x) \defto 1+\sum_{k=1}^\infty \sum_{(\u,\v)\in \tilde A^*_k(x)}\frac{1}{\sigma(\u,\v)}\prod_{m=1}^k\frac{\W^\sigma_{v_m}} {\W^\sigma_{u_m}},
$$
where
$$
A^*_k(x)=\{(\u,\v):\; x\leq u_1<v_1<\ldots <v_k<N\}
$$
and
$$
\tilde A^*_k(x)=\{(\u,\v):\; 0\leq v_1<u_1<\ldots <u_k<x\}.
$$
Then the optimal payoff to the discrete version of Problem \ref{prob3} is given by
$$
\sup_{s\in \M^{\{0,\ldots y\}}}\E_0 \lr{\prod_{t=0}^{S-1}\r_{X^{s_0}_t}}=\hat\psi (y)\defto \lrc{\sqrt{\G(y)\tilde\G^*(y)}\sum_n\A_{2n}(y)+\sqrt{\dG(y)(\Delta\tilde\G^*)(y)}\sum_n\A_{2n+1}(y)}^2,
$$
where
$$
\A_k(y)\defto \sum_{y\leq x_1<\ldots<x_k<N}\frac{1}{\sigma(\x)},
$$
$$
\dG(y)\defto \sum_{n=1}^\infty\sum_{0\leq u_1<v_1<\ldots <v_{n-1}<u_n<y}\frac{1}{\sigma(\u,\v)}\frac{\W^\sigma_{\v}}{\W^\sigma_\u}.
$$

and
$$
\Delta\tilde\G^*(y)\defto \sum_{n=1}^\infty\sum_{0\leq u_1<v_1<\ldots <v_{n-1}<u_n<y}\frac{1}{\sigma(\u,\v)}\frac{\W^\sigma_{\u}}{\W^\sigma_\v}.
$$

\endthm

\thm\label{discrete4}
The Bellman process for the dynamic version of Problem \ref{prob3} is given by
$$
V^s_t=\begin{cases}\Big (\prod\limits_{u=0}^{t-1}\r_{X^s_u}\Big )\G(X^s_t)\hat\psi(M^s_t):&\text{if }M^s_t<N,\\
\Big (\prod\limits_{u=0}^{t-1}\r_{X^s_u} \Big )\tilde\G(X^s_t):&\text{if }M^s_t=N.
\end{cases}
$$
\endthm
\subsection{Continuous-time and discrete-time with waiting}
Now we consider the cases where the  birth and death process may wait in a state and  where it forms a continuous-time Markov chain.

By solving the problem in the generality of Theorems  \ref{discrete1} to \ref{discrete4} we are are able to deal with these two cases very easily.
First, in the discrete-time case with waiting, where
$$
P_{n,n-1}=q_n;\; P_{n,n}=\e_n;\;\text{ and }P_{n,n+1}=p_n,
$$
(we stress that we take the holding probabilities $\e_n$ to be fixed and not controllable)
 we can condition on the first exit time from each state ---so that we replace 
 $P$ by $P^*$ given by
 $$
 P^*_{n,n-1}=q^*_n\defto \frac{q_n}{1-\e_n};\;  P^*_{n,n}=0;\;\text{ and }P^*_{n,n+1}=p^*_n\defto \frac{q_n}{1-\e_n}.
$$
Of course we must now modify the performance functional to allow for the time spent waiting in a state.
Thus for the additive case we must replace $f$ by $f^*$ given by 
$$f^*(n)=\frac{f(n)}{1-\e_n}.
$$
whilst in the multiplicative case we replace
$r$ by $r^*$ given by
$$
r^*(n)\defto \sum_{t=0}^\infty \e_n^t(1-\e_n)r(n)^{t+1}=\frac{(1-\e_n)r(n)}{1-\e_nr(n)}.
$$
Then in the case of a continuous-time birth and death process with birth and death rates of $\la_n$ and $\mu_n$, we obtain $P$ as the transition matrix for the corresponding jump chain --- so $P_{n,n-1}=\frac{\mu_n}{\la_n+\mu_n}$ and $P_{n,n+1}=\frac{\la_n}{\la_n+\mu_n}$ (see \cite{Feller} or \cite {R+W}). We allow for the exponential holding times by setting
$$
f^*(n)=\frac{f(n)}{\la_n+\mu_n},
$$
and
$$
r^*(n)=\frac{\la_n+\mu_n}{\a(n)+\la_n+\mu_n}
$$
Thus our results are still given by Theorems \ref{discrete1}--\ref{discrete4}

\section{Examples and some concluding remarks}\label{S-Ex}
We first consider the original time-minimisation problem with general $\sigma$.
\ex
Suppose that $f=1$ and we seek to solve Problem \ref{prob2}.
Thus $C=\emptyset$ $y=0$ and the optimal choice of $s\rq{}$ according to Theorem \ref{thm3} is
$$
s\rq{}=\frac{\sqrt 2}{\sigma}.
$$
Notice that it follows that (with this choice of scale function)
$$
\diff s(X^s_t)=\sqrt 2\diff B_t,
$$
and
$$
E_0[\ct]=s(1)^2.
$$
\endex
\ex
If we extend the previous example by assuming that $s$ is given on $[0,y)$, then we will still have $s\rq{}$ proportional to $\frac{1}{\sigma}$ on $[y,1]$ and so on this interval $s(X^s)$will behave like a multiple of Brownian Motion with partial reflection at $y$ (at least if $s\rq{}(y-)$ exists).
\endex
\ex
We now consider the additive functional case with general $f$.
Then from Theorem \ref{thm3}, the optimal choice of $s\rq{}$ is $\frac{\sqrt{2f}}{\sigma}$. With this choice of $s$, we see that
$$<s(X^s)>_t=\int_0^tf(X^s_u)\diff u,
$$
so that
$$
\E_0 \lr{\int_0^\ct f(X^s_u)du}=\E[<s(X^s)>_{\ct}].
$$
\endex
\ex
If we turn now to the discounted case and take $\alpha$ constant and $\sigma=1$, we see that the optimal choice of $s\rq{}$ is constant, corresponding to zero drift.
Thus we obtain the same optimal control for each $\alpha$. This suggests that possibly, the optimum is actually a stochastic minimum for the commute  time.
Whilst we cannot contradict this for initial position 0, the corresponding statement for a general starting position is false. 

To see this let $s_0$ correspond to drift 1 on $[0,y]$. Then a simple calculation shows that the optimal choice of $s\rq{}$ on $[y,1]$ is $\sqrt{2\a} \sqrt{ \frac{ \cosh(\sqrt{1+2\a}y) +\frac{1}{\sqrt{1+2\a}}\sinh(\sqrt{1+2\a}y)}{\cosh(\sqrt{1+2\a}y) -\frac{1}{\sqrt{1+2\a}}\sinh(\sqrt{1+2\a}y)} }
$. It is clear that this choice depends on $\a$ and hence there cannot be a stochastic minimum since, were one to exist, it would achieve the minimum in each discounted problem.
\endex
\rem
For cases where a stochastic minimum is attained in a control problem see, for example, \cite {J1} or \cite{J2}.
\endrem

\appendix

\section{Appendix: proofs}\label{Pfs}

We require the following auxiliary result.
\lem \label{Lemma1}
Let
$$B(y)\defto \int_0^y \a(u)\diff m(u)\text{ and }\tilde B(y)\defto \int_y^1 \a(u) \diff m(u),
$$
then
\eqn\label{disc2}
I_n(x)\leq \frac{(s(x)B(x))^n}{(n!)^2}\text{ and } \tilde I_n(x)\leq \frac{(\tilde s(x)\tilde B(x))^n}{(n!)^2},
\endeqn
where
$$
\tilde  s(x)\defto s(1)-s(x).
$$
\endlem
\pf

We establish the first inequality in (\ref{disc2}) by induction. The initial inequality is trivially satisfied. It is obvious from the definition that
$$
I_{n+1}(x)=\int_{v=0}^x\int_{u=0}^v \a(u)I_n(u)\diff m(u)\diff s(v),
$$
and so, assuming that $I_n(\cdot)\leq \frac{(s(\cdot)B(\cdot))^n}{(n!)^2}$:
\eqnn
I_{n+1}(x)&\leq&\int\limits_{v=0}^x\int\limits_{u=0}^v \a(u)\frac{s(u)^nB(u)^n}{(n!)^2} \diff m(u) \diff s(v)\\
&\leq& \int\limits_{v=0}^x\int\limits_{u=0}^v \a(u)\frac{B(u)^n}{(n!)^2} \diff m(u)s(v)^n \diff s(v) \text{ (since }s \text{ is increasing) }\nonumber\\
&=&\int\limits_{v=0}^x \frac{B(v)^{n+1}}{n!(n+1)!}s(v)^n \diff s(v)\nonumber\\
&\leq&\frac{B(x)^{n+1}}{n!(n+1)!}\int\limits_{v=0}^x s(v)^n \diff s(v) \text{ (since $B$ is increasing) }\nonumber\\
&=&\frac{(s(x)B(x))^{n+1}}{((n+1)!)^2},\nonumber
\endeqnn
establishing the inductive step.
A similar argument establishes the second inequality  in (\ref{disc2}).
\endpf

\pf (\textbf{of Theorem \ref{disc1}})
\begin{itemize}
\item[(i)]
Note first that $s(1)<\infty$ follows from regularity. 

We consider the case where $x\leq y$.
Now suppose that $\a$ is bounded. It follows that $\psi_y>0$ for each $y$ since $\int_0^{T_y}\a(X_u) \diff u$ is a.s. finite for bounded $\a$.
Now, setting $N_t=\exp\lrc{-\int_0^{t\wedge T_y}\a(X_u) \diff u}\psi_y(X_{t\wedge T_y})$, it is clear that $$N_t=E\lr{\exp\lrc{-\int_0^{t\wedge T_y}\a(X_u)\diff u}\Big|\F_{t\wedge T_y}}$$ and is thus a continuous martingale. Then, writing
$$
\psi_y(X_{t\wedge T_y})=\exp \lrc{\int_0^{t\wedge T_y}\a(X_u)\diff u}N_t
$$
it follows that
$$
\psi_y(X_{t\wedge T_y})-\int_0^{t\wedge T_y}\a\psi(X_u)\diff u=\int_0^t \exp\lrc{\int_0^{u\wedge T_y}\a(X_r) \diff r} \diff N_u,
$$
and hence is a martingale.
Thus we conclude that $\psi_y$ is in the domain of $\gen^y$, the extended or martingale generator for the stopped diffusion $X^{T_y}$, and 
$$
\gen^y\psi_y=\a\psi_y.
$$
Since the speed and scale measures for $X$ and $X^{T_y}$ coincide on $[0,y]$ and using the fact that $\psi_y\rq{}(0)=0$, we conclude from Theorem VII.3.12 of \cite{RY} that
\eqn\label{gen}
\psi_y(x)=\psi_y(0)+\int_{v=0}^x\int_{u=0}^vs\rq{}(v)\a(u)\psi_y(u)m\rq{}(u)\diff u \diff v\,\,\,\,\text{ for }\,x< y.
\endeqn
A similar argument establishes that
\eqn\label{gen2}
\psi_y(x)=\psi_y(1)+\int_{v=x}^1\int_{u=v}^1 s\rq{}(v)\a(u)\psi_y(u)m\rq{}(u) \diff u \diff v\,\,\,\,\text{ for }\,x> y.
\endeqn
Now either 
$$
\min(\psi_1(0),\psi_0(1))=0,
$$
in which case
$$
\E_0 \lr{\exp\lrc{-\int_0^{\ct}\a (X_t) \diff t}}=\psi_1(0)\psi_0(1)=0,
$$
or
\eqn\label{non}
\min(\psi_1(0),\psi_0(1))=c>0.
\endeqn
Suppose that (\ref{non}) holds, then (since $\psi_1$ is increasing) it follows from (\ref{gen}) that
\eqnn
\psi_1(1-)&\geq& c+\int_{v=0}^1\int_{u=0}^vs\rq{}(v)c\a(u)m\rq{}(u) \diff u \diff v\\
&=&c\lrc{1+\int_{u=0}^1\int_{v=u}^1s\rq{}(v)\a(u)m\rq{}(u) \diff u \diff v}\nonumber\\
&\geq& c\lrc{1+\lr{ s(1)-s\lrc{\half}}\int_0^{\half}\a(u)m\rq{}(u) \diff u}.\nonumber
\endeqnn

Similarly, we deduce that
$$
\psi_0(0+)\geq c\lrc{1+s\lrc{\half}\int_{\half}^1\a(u)m\rq{}(u)\diff u}.
$$
Thus, if (\ref{bound}) fails, (\ref{non}) cannot hold (since if (\ref{bound}) fails then at least one of $\int_0^{\half}\a(u)m\rq{}(u)\diff u)$ and $\int_{\half}^1\a(u)m\rq{}(u) \diff u)$ is infinite) and so we must have $\psi_1(0)\psi_0(1)=0$.

To deal with unbounded $\a$, take a monotone, positive  sequence $\a_n$ increasing to $\a$ and take limits.

\item[(ii)]Suppose now that (\ref{bound}) holds.
Setting
$$
G(x)=\frac{\psi_1(x)}{\psi_1(0)},
$$
we see that
$G$ satisfies equation (\ref{gen}) with $G(0)=1$.

Convergence of the series $\sum I_n$ and $\sum \tilde I_n$ follows from the bounds on $I_n$ and $\tilde I_n$ given in  Lemma \ref{Lemma1}

Now by iterating equation (\ref{gen}) we obtain 
\eqnno
G(x)&=&\sum_{k=0}^{n-1} I_k(x)\\
&+&\int\limits_{0\leq u_1\leq v_1\leq u_2\ldots \leq v_n\leq x}\a(u_1)\ldots \a(u_n)G(u_n)\diff m(u_1)\ldots \diff m(u_n)\diff s(v_1)\ldots \diff s(v_n).
\endeqnno
Since $G$ is bounded by $\frac{1}{\psi_1(0)}$ we see that
$$
0\leq G(x)-\sum_{k=0}^{n-1} I_k(x)\leq \frac{1}{\psi_1(0)}I_n(x).
$$
A similar argument establishes that
$$
0\leq \tilde G(x)-\sum_{k=0}^{n-1} \tilde I_k(x)\leq \frac{1}{\psi_0(1)}\tilde I_n(x).
$$
and so we obtain (\ref{disc}) by taking limits as $n\rightarrow \infty$.\endpf
\end{itemize}

\pf (\textbf{of Theorem \ref{thm3}})
Note first that, from Theorem \ref{time},
\eqnn\nonumber
\E_0\lr{\int_0^\ct f(X^s_t)\diff t}&=&\phi_1(0)+\phi_0(1)=\int\limits_{v=0}^1\int\limits_{u=0}^1f(u)s(\diff v)m(\diff u)\\
=s_0(C)I^{s_0}(C)&+&\int\limits_{C^c}[I^{s_0}(C)s(\diff v)+s(C)f(v)m( \diff v)]\nonumber\\
&&+\half \int\limits_{C^c}\int\limits_{C^c}[f(u)s(\diff v)m(\diff u)+f(v)m(\diff v)s(\diff u)],\label{add2}
\endeqnn
where the factor $\half$ in the last term in (\ref{add2})  arises from the fact that we have symmetrised the integrand.
Now, for $s\in \M^C_{s_0}$, we can rewrite (\ref{add2}) as
\eqnn \nonumber
\E_0\lr{\int_0^\ct f(X_t)dt}=s_0(C)I^{s_0}(C)&+\,\,\,\int\limits_{C^c}\lr{I^{s_0}(C)s\rq{}(v)+s_0(C)\frac{2f(v)}{\sigma^2(v)s\rq{}(v)}} \diff v\\
&+\,\,\,\half\int\limits_{C^c}\int\limits_{C^c}\lr{\frac{2f(u)}{\sigma^2(u)}\frac{s\rq{}(v)}{s\rq{}(u)}+\frac{2f(v)}{\sigma^2(v)}\frac{s\rq{}(u)}{s\rq{}(v)}} \diff u \diff v.\label{add3}
\endeqnn
We now utilise the very elementary fact that for 
$a,b\geq 0$,
\eqn\label{ineq}
\inf_{x>0}\lr{ax+\frac{b}{x}}=2ab\text{ and if }a,b>0\text{ this is attained at }x=\sqrt{\frac{b}{a}}.
\endeqn

Applying this to the third term on the right-hand-side of (\ref{add3}), we see from (\ref{ineq}) that it is bounded below by $\int\limits_{C^c}\int\limits_{C^c}\sqrt{\frac{4f(u)f(v)}{\sigma^2(u)\sigma^2(v)}}\diff u \diff v=\lrc{\int_{C^c}\sqrt{\frac{2f(u)}{\sigma^2(u)}} \diff u}^2=J^2(C^c)$ and this bound is attained when $s\rq{}(x)$ is a constant multiple of $\sqrt{\frac{2f(x)}{\sigma^2(x)}}$ a.e. on $C^c$. 

Turning to the second term in  (\ref{add3}) we see from (\ref{ineq}) that it is bounded below by $\int_{C^c}2\sqrt{s_0(C)I^{s_0}(C)}\sqrt{\frac{2f(v)}{\sigma^2(v)}} \diff v=2\sqrt{s_0(C)I^{s_0}(C)}J(C^c)$ and this is attained when $s\rq{}(x)=\sqrt{\frac{2f(x)}{\sigma^2(x)}}\sqrt{\frac{s_0(C)}{I^{s_0}(C)}}$ a.e. on $C^c$. 

Thus, we see that the infimum of the RHS of (\ref{add3}) is attained by setting $s\rq{}(x)$ equal to $\sqrt{\frac{s_0(C)}{I^{s_0}(C)}}\sqrt{\frac{2f(x)}{\sigma^2(x)}}$ on $C^c$ and this gives the stated value for the infimum.
\endpf

\pf (\textbf{of Theorem \ref{discstat}})
\begin{itemize}
\item[(i)]This is proved in the same way as equation (\ref{disc}) in Theorem \ref{disc1}.
\item[(ii)] First we define 
\eqn
I_n^*(x)\defto \int\limits_{x\leq u_1\leq v_1\leq\ldots v_n\leq 1}
\a(v_1)\ldots \a(v_n) \diff s(u_1)\ldots  \diff s(u_n) \diff m(v_1)\ldots  \diff m(v_n);
\endeqn
and
\eqn
\G^*(x)\defto \sum_{n=0}^\infty I_n^*(x).
\endeqn
To prove (ii) we use the following representations (which the reader may easily verify):
\eqn
I_n(1)=\sum_{m=0}^nI_m(y) I^*_{n-m}(y)-\tilde \sigma^2(y)\sum_{m=1}^n I_m\rq{}(y)(I^*_{n-m})\rq{}(y),
\endeqn
and
\eqn
\tilde I_n(0)=\sum_{m=0}^n\tilde I_m(y)\tilde I^*_{n-m}(y)-\tilde \sigma^2(y)\sum_{m=1}^n \tilde I_m\rq{}(y)(\tilde I^*_{n-m})\rq{}(y)
\endeqn

It follows from these equations that
\eqnn\label{add4}
&\E_0& \lr{ \exp\lrc{- \int_0^\ct \a(X_t) \diff t}}\\
&=&\bigl[G(y)G^*(y)-\tilde \sigma^2(y)G\rq{}(y)(G^*)\rq{}(y)\bigr]\bigl[\tilde G(y)\tilde G^*(y)-\tilde \sigma^2(y)\tilde G\rq{}(y)(\tilde G^*)\rq{}(y)\bigr].\nonumber
\endeqnn
Now essentially the same argument as in the proof of Theorem \ref{thm3} will work as follows.
Multiplying out the expression on the RHS of (\ref{add4}) we obtain the sum of the three terms:
\begin{itemize}
\item[(a)]$\half G(y)\tilde G^*(y)\sum\limits_{m\geq 0 ,n\geq 0}[\tilde I_n(y)I^*_m(y)+\tilde I_m(y)I^*_n(y)]$
\item[(b)]$\half G\rq{}(y)\tilde (G^*)\rq{}(y)\sum\limits_{m\geq 0 ,n\geq 0}[\tilde I\rq{}_n(y)(I^*_m)\rq{}(y)+\tilde I\rq{}_m(y)(I^*_n)\rq{}(y)]$ ;
and
\item[(c)]$\sum\limits_{m\geq 1,n\geq 0}[G(y)(\tilde G^*)\rq{}(y)I^*_n(y)\tilde I_m\rq{}(y)+G\rq{}(y)\tilde G^*(y)(I^*_m)\rq{}(y)\tilde I_n(y)]$,
\end{itemize}
where in the first two terms we have symmetrised the sums.

Using (\ref{disc}), the sum in (c) becomes
\eqnn\label{ceq}
\sum_{m\geq 1,n\geq 0}&\int\limits_{D_{m,n}(y)}& \left [G(y)(\tilde G^*)\rq{}(y)\frac{t\rq{}(v_1)\ldots t\rq{}(v_n)t\rq{}(w_1)\ldots t\rq{}(w_m)}{t\rq{}(u_1)\ldots t\rq{}(u_n)t\rq{}(z_1)\ldots t\rq{}(z_{m-1})} \right.\\
&&+ \left. G\rq{}(y)\tilde G^*(y)\frac{t\rq{}(u_1)\ldots t\rq{}(u_n)t\rq{}(z_1)\ldots t\rq{}(z_{m-1})}{t\rq{}(v_1)\ldots t\rq{}(v_n)t\rq{}(w_1)\ldots t\rq{}(w_m)}\right ]\diff \tilde\la (\u,\v,\underline w,\z)\nonumber 
\endeqnn
where
\eqnno
D_{m,n}(x)=\{(\underline u, \underline v,\underline w,\underline z)\in \R^n\times\R^n\times\R^{m}\times\R^{m-1}&:&\; x\leq u_1\leq v_1\leq\ldots v_n\leq 1;\\
\text{ and }&&x \leq w_1\leq z_1\leq\ldots w_m\leq 1\},
\endeqnno 
$t$ is the measure with Radon-Nikodym derivative $t\rq{}=\tilde \sigma s\rq{}$,
and $\tilde \la$ denotes the measure with Radon-Nikodym derivative $\frac{1}{\tilde\sigma}$.
Clearly each term in the sum in (\ref{ceq}) is minimised by taking $t\rq{}$ constant and equal to $\sqrt{\frac{G(y)\H \rq{}(y)}{G\rq{}(y)\H (y)}}$ a.e. on $[y,1]$.

The first two terms, (a) and (b), are each minimised by taking $t\rq{}$ constant a.e. on $[y,1]$. Substituting this value for $t\rq{}$ back in we obtain the result.
\end{itemize}
\endpf

\pf (\textbf{of Theorem \ref{dsc1}})
Consider the candidate Bellman process $v_t$. 
Using the fact that
\eqn
N_t\defto\int_0^{t\wedge T_1}f(X^{s_0}_u)\diff u-\phi_{X^{s_0}_{t\wedge T_1}}(0)=\E\lr{\int_0^{T_1}f(X^{s_0}_u)\diff u|\F_{t\wedge T_1}}
\endeqn
is a martingale,
\eqn
N_t\rq{}\defto\phi_0({X^{s_0}_t})+\int_0^t f(X^{s_0}_u)\diff u
\endeqn
is equal to $\E\lr{\int_0^\ct f(X^{s_0}_u) \diff u\big |\F_t}$ on the stochastic interval $[[T_1,\ct]]$, and hence is a martingale on that interval, $M^{s_0}$ is a continuous, increasing process, and $\phi_1(0)+\phi_0(1)=2s_0(1)I^{s_0}(1)$ (so that $v$ is continuous at $\T_1(X^{s_0})$):

\eqnno
\diff v_t&=&4\biggl(\sqrt{{s_0}(M^{s_0}_t)I^{s_0}(M^{s_0}_t)}+J(M^{s_0}_t)\biggr)\\
&&\times\biggl[\half s_0\rq{}(M^{s_0}_t)\sqrt{\frac{I^{s_0}(M^{s_0}_t)}{{s_0}(M^{s_0}_t)}}
+\half (I^{s_0})\rq{}(M^{s_0}_t)\sqrt{\frac{{s_0}(M^{s_0}_t)}{I^{s_0}(M^{s_0}_t)}}-\sqrt{\frac{2f(M^{s_0}_t)}{\sigma^2(M^{s_0}_t)}}\biggr] \diff M^{s_0}_t\\
&&+\diff N_t1_{(M^{s_0}_t<1)}+\diff N_t\rq{}1_{(M^{s_0}_t=1)}\\
&=&4\biggl(\sqrt{{s_0}(M^{s_0}_t)I^{s_0}(M^{s_0}_t)}+J(M^{s_0}_t)\biggr)\\
&\times&\biggl[\half\biggl(s\rq{}(M^{s_0}_t)\sqrt{\frac{I^{s_0}(M^{s_0}_t)}{{s_0}(M^{s_0}_t)}}
+\frac{1}{s\rq{}(M^{s_0}_t)}\frac{2f(M^{s_0}_t)}{\sigma^2(M^{s_0}_t)}\sqrt{\frac{{s_0}(M^{s_0}_t)}{I^{s_0}(M^{s_0}_t)}}\biggr)-\sqrt{\frac{2f(M^{s_0}_t)}
{\sigma^2(M^{s_0}_t)}}\biggr]\diff M^{s_0}_t\\
&+&\diff N_t1_{(M^{s_0}_t<1)}+\diff N_t\rq{}1_{(M^{s_0}_t=1)}
\endeqnno
Now, since $s_0$, $I^{s_0}$ and $J$ are non-negative it follows from (\ref{ineq}) that
$$\diff v_t\geq \diff \bar N_t,
$$
where
$$
\diff \bar N_t=\diff N_t1_{(M^{s_0}_t<1)}+\diff N_t\rq{} 1_{(M^{s_0}_t=1)},
$$
with equality if 
\eqn\label{opts}
s\rq{}(M^{s_0}_t)=\sqrt{\frac{2s_0(M^{s_0}_t)f(M^{s_0}_t)}{I^{s_0}(M^{s_0}_t)\sigma^2(M^{s_0}_t} }.
\endeqn
Then the usual submartingale argument (see, for example \cite{Oks} Chapter 11), together with the fact that $v$ is bounded by assumption (\ref{ass1})) gives us 
(\ref{soln1}). 

It is easy to check that $s$ given by (\ref{opts}) is in $\M^{s_0}_{M^{s_0}_t}$.
The fact that the optimal choice of $s$ satisfies (\ref{soln2}) follows on substituting $s\rq{}(x)=\sqrt{\frac{s_0(y)f(x)}{I^{s_0}(y)\sigma^2(x)} }$ in the formulae for $s$ and $I^{s}$ and observing that the ratio $\frac{s_0(x)}{I^{s_0}(x)}$ is then constant on $[y,1]$.
\endpf

\pf (\textbf{of Theorem \ref{dsc2}})
The proof is very similar to that of Theorem \ref{dsc1}. Note that $\psi$ is continuous at the point $(1,1)$.

Thus, for a suitable bounded martingale $n$,
\eqnno
\diff \pay_t&=&\exp\lrc{-\int_0^t \alpha(X^{s_0}_u)\diff u}\psi_y(X^{s_0}_t,M^{s_0}_t)\diff M^{s_0}_t1_{(M^{s_0}_t<1)}+\diff n_t\\
&=&\exp\lrc{-\int_0^t \alpha(X^{s_0}_u)\diff u}G(X^{s_0}_t)\psi\rq{}(M^{s_0}_t)1_{(M^{s_0}_t<1)}+\diff n_t\\
&=&-2\exp\lrc{-\int_0^t \alpha(X^{s_0}_u)\diff u}G(X^{s_0}_t)\lrc{\sqrt{G{\H }}\cosh \A(M^{s_0}_t)+\sqrt{\tilde\sigma^2G\rq{}\H \rq{}}\sinh \A(M^{s_0}_t)}^{-3}\\
&&\times\lr{ \lr{ \lrc{\sqrt{G\H }}\rq{}-\sqrt{G\rq{}\H \rq{}}}\cosh \A(M^{s_0}_t)+\lrc{\sqrt{\tilde\sigma^2G\rq{}\H \rq{}}}\rq{}-\sqrt{\tilde\sigma^2G\H }\sinh \A(M^{s_0}_t)}+\diff n_t.
\endeqnno
Now
$$
\lrc{\sqrt{G\H }}\rq{}=\half \H \rq{}\sqrt{\frac{G}{\H }}+\half G\rq{}\sqrt{\frac{\H }{G}}\geq \sqrt{G\rq{}\H \rq{}}\text{ using (\ref{ineq})},
$$
with equality attained when 
\eqn\label{d4}
\sqrt{G\rq{}\H \rq{}}=\sqrt{G\H }.
\endeqn
Similarly, defining $m^{\a}$ by setting $\diff m^{\a}=\frac{\diff m}{\a}$,
\eqnno
\lrc{\sqrt{\tilde\sigma^2G\rq{}\H \rq{}}}\rq{}&=&\lrc{\sqrt{(\frac{G\rq{}}{s\rq{}})(\tilde\sigma^2\H \rq{}}}\rq{}=\lrc{\sqrt{\frac{\diff G}{\diff s}\frac{\diff \H }{\diff \am }}}\rq{}\\
&=&\half (\am)\rq{}\frac{\diff^2G}{\diff \am ds}\sqrt{\frac{\frac{\diff \H }{\diff \am}}{\frac{\diff G}{\diff s}}}+\half s\rq{}\frac{\diff^2\H }{\diff s \diff \am }\sqrt{\frac{\frac{\diff G}{\diff s}}{\frac{\diff \H }{\diff \am }}}\\
&=&\half \lrc{\frac{1}{\tilde\sigma^2s\rq{}}\frac{\diff^2G}{\diff \am \diff s}\sqrt{\frac{\frac{\diff \H }{\diff \am }}{\frac{\diff G}{\diff s}}}+ s\rq{}\frac{\diff^2\H }{\diff s \diff \am}\sqrt{\frac{\frac{\diff G}{\diff s}}{\frac{\diff \H }{\diff \am }}}}\\
&=&\half \lrc{\frac{1}{\tilde\sigma^2s\rq{}}G\sqrt{\frac{\frac{\diff \H }{\diff \am }}{\frac{\diff G}{\diff s}}}+s\rq{}\H \sqrt{\frac{\frac{\diff G}{\diff s}}{\frac{\diff \H }{\diff \am }}}}\\
&\geq&\sqrt{\frac{G\H } {\tilde\sigma^2} }
\endeqnno
with equality when
\eqn\label{d5}
s\rq{}=\frac{1}{\tilde\sigma}\sqrt{\frac{G{\frac{\diff \H  }{\diff \am }}} {\H  \frac{\diff G}{\diff s}}   }
\endeqn
Now we can easily see (by writing $\frac{dG}{ds}=\frac{1}{s\rq{}}G\rq{}$ and $\frac{d\H  }{d\am}=\tilde\sigma^2s\rq{}\H  \rq{}$) that (\ref{d5}) implies (\ref{d4}) so the standard supermartingale argument establishes that
$$
V_t=\pay_t.
$$
That the optimal choice of $s\rq{}$ is as given in (\ref{d3}) follows on observing that, with this choice of $s\rq{}$,
$$
(\tilde\sigma(G\rq{}\H  -G\H  \rq{}))\rq{}(x)=0\text { for }x\geq y,
$$
and
$$
G\rq{}(y)\H  (y)-G(y)\H  \rq{}(y)=0.
$$
\endpf


\begin{thebibliography}{99}

\bibitem{Ass}Assing, S and Schmidt, W M: \lq\lq Continuous Strong Markov Processes in Dimension One'', Springer, Berlin-Heidelberg-New York (1998).
\bibitem{ARR}Atchade, Y, Roberts, G O and J. S. Rosenthal:  Towards Optimal Scaling of Metropolis-Coupled Markov Chain Monte Carlo', {\it Statistics and Computing}, 21, 4, 555-568, 2011.
\bibitem{Barlow} Barlow, M T\lq\lq Random Walks and Heat Kernels on Graphs'', Cambridge University Press, Cambridge, 2017.
\bibitem{J2}Connor, S and Jacka, S: Optimal co-adapted coupling for the symmetric random walk on the hypercube, {\it J. Appl. Probab.} 45, 703-713, 2008.
\bibitem{CRRST} Chandra, A K,  Raghavan, P, Ruzzo, V L, Smolensky, R, and Tiwari, P.: The electrical
resistance of a graph captures its commute and cover times. {\it Proceedings of the 21st
ACM Symposium on theory of computing}, 1989.
\bibitem{E2}Englebert, H J: Existence and non-existence of solutions of one-dimensional stochastic equations, {\it Probability and Mathematical Statistics},Vol. 20 (2), 343 - 358, 2000. 
\bibitem{E}Englebert, H J and Schmidt, W, On solutions of one-dimensional stochastic differential equations without drift, {\it Z. Wahrsch. Verw. Gebiete} 68, 287 - 314, 1985.
\bibitem{Feller}Feller, W: \lq\lq An Introduction To Probability Theory and
its Application Vol II''. 2nd Edition. Wiley New York, 1971.
\bibitem{FSoner} Fleming, W. H. and Soner H. M.: \lq\lq Controlled Markov Processes and Viscosity Solutions'',  Springer-Verlag, Berlin, 1993.
\bibitem{IM}Ito, K and McKean, H P: \lq\lq Diffusion Processes and their Sample Paths\rq\rq{}, Springer, Berlin-Heidelberg-New York, 1974.
\bibitem{J1}Jacka, S:, Keeping a satellite aloft: two finite fuel stochastic control models, {\it J. Appl. Probab.} 36, 1-20, (1999).
\bibitem{Oks}Oksendal, B: \lq\lq{}Stochastic Differential Equations\rq\rq{}, 6th Edition, Springer, Berlin-Heidelberg-New York, 2003.
\bibitem{RY}Revuz, D and Yor, M: \lq\lq Continuous martingales and Brownian motion'', 3rd Edn. Corrected 3rd printing.
Springer, Berlin-Heidelberg-New York, 2004.
\bibitem{R+W}Williams, D: \lq\lq Diffusions, Markov processes, and martingales: Vol. I''. 
Wiley, New York, 1979.









\end{thebibliography}
\end{document}